\documentclass{article}
\usepackage[pdftex]{hyperref}
\usepackage{tikz-cd, url}

\hypersetup{
  colorlinks   = true, 
  urlcolor     = gray, %Color for external hyperlinks
  linkcolor    = red, %Color of internal links
  citecolor   =  blue%Color of citations
}
\usepackage{amsmath,amsfonts,amsthm,amssymb,graphicx,xcolor}
\usepackage[all,2cell,ps]{xy}

\usepackage{todonotes}
\theoremstyle{plain}
\newtheorem{thm}{Theorem}[section]

\newtheorem{prop}[thm]{Proposition}
\newtheorem{cor}[thm]{Corollary}
\newtheorem{conj}[thm]{Conjecture}

\theoremstyle{definition}

\newtheorem{ex}[thm]{Example}

\newtheorem{rem}[thm]{Remark}

\newcommand{\Z}{\mathbb Z}
\newcommand{\R}{\mathbb R}

\newcommand{\C}{\mathbb C}

\newcommand\OO{{\mathcal O}}

\newcommand\inv{^{-1}}

\DeclareMathOperator{\Pic}{Pic}

\DeclareMathOperator{\spec}{Spec}

\usepackage{mathtools}

\DeclarePairedDelimiter\abs{\lvert}{\rvert}
\DeclarePairedDelimiter\norm{\lVert}{\rVert}

\makeatletter
\let\oldabs\abs
\def\abs{\@ifstar{\oldabs}{\oldabs*}}
\let\oldnorm\norm
\def\norm{\@ifstar{\oldnorm}{\oldnorm*}}
\makeatother
\newenvironment{pf}{\begin{proof}}{\end{proof}}

\title{On the Hopf Problem and a Conjecture of Liu-Maxim-Wang}
 \author{\small{Luca F. Di Cerbo}\footnote{Supported in part by NSF grant DMS-2104662} \\ \scriptsize{University of Florida} \\ \footnotesize{\textsf{ldicerbo@ufl.edu}} \and 
\small{Rita Pardini} \footnote{Partially supported by PRIN  2017SSNZAW\_004. Member of Gnsaga-INdAM.}  \\ 
\scriptsize{Università di Pisa}\\ \footnotesize{\textsf{rita.pardini@unipi.it}}}
\date{}
	
\begin{document}

\maketitle

\begin{abstract}
 We discuss an approach towards the Hopf problem for aspherical smooth projective varieties recently proposed by Liu, Maxim, and Wang in \cite{LMW17b}. In complex dimension two, we point out that this circle of ideas suggests an intriguing conjecture regarding the geography of aspherical surfaces of general type.
\end{abstract}
%%%%%%%%%%%%%%%%%%%%
\vspace{10cm}
\tableofcontents\quad\\

\vspace{1cm}

%%%%%%%%%%%%%%%%%%%%%%%%%%%%%%%%
\section{Introduction}
%%%%%%%%%%%%%%%%%%%%%%%%%%%%%%%%

A long-standing and important problem in geometry  is a conjecture of Hopf on the sign of the Euler characteristic of aspherical manifolds. 

\begin{conj}[Hopf Conjecture]\label{Hopf}
	If $X$ is a closed aspherical manifold of real dimension $2n$, then:
	\begin{equation*}
	(-1)^n\chi_{\rm top}(X)\geq 0.
	\end{equation*}
	
\end{conj}

For more details about this problem, we refer to M. Berger's panoramic book on Riemannian geometry and to S.-T. Yau's authoritative list of problems in geometry, see in particular \cite[Chapter 12, Note 12.3.1.1]{Berger} and \cite[Section VII, Problem 10]{SchoenY}. Conjecture \ref{Hopf} is true when $n=1$ thanks to the uniformization theorem for Riemann surfaces. Interestingly, this problem is still open when $n=2$. Moreover, Conjecture \ref{Hopf} is still wide open even in the realm of aspherical, smooth, projective varieties! This is particularly surprising as, in the projective case, one can use a large variety of tools coming from algebraic geometry, on top of the usual differential geometry, geometric analysis, and geometric topology approaches to such conjecture.

The literature around the Hopf problem is very vast and it continues to grow, we refer the interested reader to \cite{Chern}, \cite{DX84}, \cite{Gro}, \cite{Ko}, \cite{JZ}, \cite{PS13}, \cite{LMW17b}, \cite{arapura}, \cite{DS22}, \cite{Maxim} for a selection of diverse contributions on this topic over a seven decades period. In an interesting recent paper \cite{LMW17b}, Y. Liu, L. Maxim, and B. Wang connect the Hopf problem for smooth projective varieties with the well-known Shafarevich conjecture in algebraic geometry. More precisely, in \cite[Theorem 1.8]{LMW17b} they observe that if the Shafarevich conjecture is true then any aspherical, smooth, projective variety must have Stein topological universal cover. They also conjecture (see \cite[Conjecture 6.3]{LMW17b}) that varieties with Stein universal cover have \emph{nef} cotangent bundle, and show that the sign of the Euler characteristic of a smooth variety with nef cotangent bundle satisfies the statement of Conjecture \ref{Hopf} (see \cite[Proposition 3.6]{LMW17b}). Therefore, they suggest a new interesting approach to Conjecture \ref{Hopf} via the Shafarevich conjecture and the study of the nefness properties of the holomorphic cotangent bundle.

In \cite{Yiyu}, Y. Wang (a graduate student of B. Wang) addresses Conjecture 6.3 in \cite{LMW17b} and produces some examples of smooth projective varieties with Stein universal topological cover, but with \emph{non-nef} holomorphic cotangent bundle. The goal of this paper is to show that such examples are \emph{not} aspherical. The question of the asphericity of such examples was raised during the recent \textbf{American Mathematical Society} Special Session on  ``Singer-Hopf Conjecture in Geometry and Topology'', that was part of the \textbf{Spring Southeastern Sectional Meeting}, Georgia Institute of Technology, Atlanta, GA, March 18-19, 2023. During this special session, L. Maxim gave a lecture on the program towards the Hopf conjecture contained in \cite{LMW17b}, \cite{Maxim}, and Y. Wang presented the examples contained in \cite{Yiyu}.

We provide several proofs of the \emph{non-asphericity} of the examples described in \cite{Yiyu}. All such proofs crucially rely upon results of Tovena and Pardini \cite{PT95} concerning the fundamental group of finite abelian branched covers. In complex dimension two, we also offer a proof entirely based on complex analysis. Indeed, we show that Y. Wang's examples are deformation equivalent to smooth projective varieties containing smooth rational curves. These spaces are then non-aspherical in a very strong sense, and they are deformation equivalent to smooth projective varieties with non-Stein universal cover.

Concluding, the approach to Conjecture \ref{Hopf} via the Shafarevich conjecture and the study of the nefness properties of the holomorphic cotangent bundle is still a very viable one. The next intriguing question to explore is whether smooth, aspherical, projective varieties have \emph{nef} holomorphic cotangent bundle (see \cite[Conjecture 6.4]{LMW17b}). Even if Conjecture \ref{Hopf} is known to be true for aspherical complex surfaces thanks to the Kodaira-Enriques classification, the nefness of their cotangent bundle is currently unknown. Indeed, in Section \ref{surfaces} we point out several intriguing questions on the geography of aspherical surfaces of general type suggested by this circle of ideas.

\noindent\textbf{Acknowledgments}. 
The first named author thanks the participants of the AMS Special session on ``Singer-Hopf Conjecture in Geometry and Topology'' at the 2023 Spring Southeastern Sectional Meeting  for a very stimulating exchange of ideas around the Hopf problem. He also thanks Yiyu Wang for a useful correspondence. The authors thank J\'anos Koll\'ar for useful bibliographical suggestions and for pertinent comments on the manuscript (see Remark \ref{rem: kollar}), and the referee for useful remarks.\\

\section{Main Result -- Several Proofs}\label{Proofs}\label{sec: main}

A smooth closed manifold $X$ is said to be \emph{aspherical} if its topological universal cover $\widetilde{X}$ is contractible. This is equivalent to the vanishing of the higher homotopy groups $\pi_i(X)$, for $i\geq 2$. Recall that for $i\geq 2$ we have the isomorphism $\pi_i(X) \cong \pi_i(\widetilde{X})$ (see \cite[Proposition 4.1]{Hat}). We refer in general to Chapter $4$ in Hatcher's book \cite{Hat} for the relevant definitions and the basic homotopy theory results.

%\ritainline{qui sotto ho editato un po' e anche cambiato la notazione per usare quella a cui sono abituata, spero non ti dispiaccia}
The examples constructed in \cite{Yiyu} are easily described. Start with an abelian variety of dimension $n$ (a projective complex torus of real dimension $2n$), say $Y$. The topological universal cover is then given by 
\[
\pi\colon \C^n\to Y, \quad Y\cong  \C^n/ \Lambda, \quad \Lambda\cong \Z^{2n}, 
\]
where $\widetilde{Y}=\C^n$ is clearly contractible and Stein. Next, one takes an  ample line bundle, say $L$, on $Y$ and a smooth divisor  $B$ in the linear system $|2L|$. One can define on $\mathcal O_Y\oplus L^{-1}$ a $\Z_2$-action by declaring the first summand to be invariant and the second one anti-invariant,  and then an   $\OO_Y$-algebra  structure compatible with this action using the map $L\inv\otimes L\inv \cong \OO_Y(-B)\hookrightarrow \OO_X$. Setting $X:=\spec(\mathcal O_Y\oplus L^{-1})$ defines a degree two finite morphism of projective varieties $f\colon X\to Y$, that we call the {\em double cover given by the relation $2L\sim B$} ($\sim$ denotes linear equivalence). The variety $X$ is smooth, since  $B$ is, and standard formulae for double covers give $K_X=f^*(K_Y\otimes L)=f^*(L)$, where $K_X$ denotes as usual the canonical bundle.  So $K_X$ is ample and  $X$ is minimal of general type.  As shown in \cite[Corollary 2]{Yiyu}, the holomorphic cotangent bundle $\Omega^1_{X}$ is not nef, while the topological universal cover $\widetilde{X}$ is Stein.

\begin{rem}\label{rem: severi-line}
The above examples have maximal Albanese dimension, namely their Albanese map is generically finite (it is easily seen to coincide with $f$, see \cite[\S 2.4 (d)]{geography}); they satisfy  equality in the generalized Severi inequality $K_X^n\ge 2n! \chi(K_X)$ for varieties of general type and maximal Albanese dimension (\cite{barja-severi}, \cite{zhang-severi}). Recall that $\chi(K_X)=(-1)^n\chi(\OO_X)$ by Serre duality and that $\chi(K_X)\ge 0$ if $X$ is of Albanese general type  by the Generic Vanishing Theorem (Corollary of \cite[Theorem~1]{GL1}). Thus, they satisfy the inequality $\chi(K_{X})>0$. Conjecturally, this inequality holds true for a general class of varieties, see for example \cite[Conjecture 18.12.1]{KBook}.
\par\noindent 
In complex dimension two,  it is proven in \cite{BPS} that the minimal surfaces $S$ of general type with maximal Albanese dimension and $K^2_S=4\chi(\OO_S)$ have $h^1(\OO_S)=2$ and  are precisely the (minimal resolutions of) double covers of their Albanese surface  branched on an ample divisor $B$ with at most A-D-E singularities (see \cite[\S~II.8]{bpv} for the definition).
The cover is smooth (and $K_S$ is ample) precisely when $B$ is smooth, otherwise it has rational double points and is the canonical model of a surface of general type.
\end{rem}

\begin{thm}\label{main}
If $n\ge 2$, the $n$-dimensional smooth projective variety $X$ is not aspherical.
\end{thm}
\begin{pf}
Assume by contradiction that $\widetilde{X}$ is contractible. Since $\widetilde{Y}=\C^n$ is contractible and Stein, by Corollary 3.4 in \cite{PT95} we know that
\[
f_{*}\colon \pi_1(X)\to\pi_1(Y)\cong \Z^{2n}
\]
is an isomorphism. By Whitehead's theorem (see \cite[Theorem 4.5]{Hat}), a map between connected CW complexes that induces isomorphisms on all homotopy groups is necessarily a homotopy equivalence. Thus, we conclude that $f\colon X\to Y$ is a homotopy equivalence. 
At this point there are several different arguments that can be used to obtain a contradiction, showing that $X$ is not aspherical:
%\ritainline{edita e formatta come preferisci}
\begin{itemize}
\item (Topological arguments)\newline
Assume by contradiction that $X$ is aspherical. $X$ and $Y$ have a choice of orientations that determines isomorphisms  $H_{2n}(X; \Z)\cong \Z$ and 
$H_{2n}(Y; \Z)\cong \Z$. Since $f\colon X\to Y$ is a double branched cover,  with these choices the induced map  $f_*\colon H_{2n}(X; \Z)\to  H_{2n}(Y; \Z) $ is multiplication by $2$. Thus, $f$ cannot  be a homotopy equivalence.  We have reached a contradiction and $X$ cannot be aspherical.

Alternatively, we could argue with the Euler characteristic instead of the degree of the map. Since $f$ is a homotopy equivalence and $Y$ is a $2n$-dimensional real torus  we have $\chi_{\rm top}(X)=\chi_{\rm top}(Y)=0$.
On the other hand, we are going to show   $(-1)^n\chi_{\rm top}(X)>0$.  First, the double  cover $f\colon X\to Y$ induces an isomorphism of the ramification divisor with the branch divisor $B$ and restricts to a topological cover of degree 2 of $Y\setminus B$. By Mayer--Vietoris for the Euler characteristic, we then have:
\[
\chi_{\rm top}(X)=2\chi_{\rm top}(Y\setminus B)+\chi_{\rm top}(B)=2\chi_{\rm top}(Y)-\chi_{\rm top}(B)=-\chi_{\rm top}(B)
\]
Since $B$ is an ample divisor in $Y$, by Lefschetz hyperplane theorem (see for example \cite[Theorem 3.1.17]{Laz}) we have
that the map in cohomology associated to the inclusion of $B$ in $Y$  
\[
H^{i}(Y; \Z)\to H^{i}(B; \Z)
\]
is an isomorphism for $i\leq n-2$ and an injection for $i=n-1$. Since $Y$ is a torus of real dimension $2n$, we have $b_{i}(Y)=\binom{2n}{i}$ for $0\leq i\leq 2n$. We conclude
\[
(-1)^n\chi_{\rm top}(X)=(-1)^{n-1}\chi_{\rm top}(B)\geq \binom{2n}{n}-\binom{2n}{n-1}>0.
\]

\item (Hodge theoretic argument)\newline
 The map $f^*\colon H^k(Y;\C)\to H^k(X; \C)$ is an isomorphism, since $f$ is a homotopy equivalence, and preserves the Hodge decomposition, since $f$ is a morphism. So $X$ and $Y$ have the same Hodge numbers. In particular $h^{n,0}(Y)=h^{n,0}(X)=1$. On the other hand, since 
$K_X=f^*(K_Y\otimes L)=f^*(L)$, the projection formula for a finite flat morphism gives
\[
f_*K_X=f_{*}f^{*}(L)=f_*(\OO_Y)\otimes L=(\OO_Y\oplus L^{-1})\otimes L=  L\oplus\OO_{Y}.
\]
Since $h^0(K_X)=h^0(f_*K_X)$, we then have
$$h^{n,0}(X)=1+h^0(L)=1+\frac{L^n}{n!}>1,$$
where the last equality follows by Riemann-Roch and Kodaira vanishing, as $L$ is ample by assumption. We have reached a contradiction.

\item (Rigidity argument)\newline
A classical rigidity result of Catanese (\cite[Proposition 4.8]{Catanese}) implies that $X$ is an abelian variety, since it is  K\"ahler with $b_1(X)=2n$  and the cohomology algebra $H^{\bullet}(X;\Z)$ is isomorphic  to the exterior algebra $\wedge^{\bullet}H^1(X;\Z)$. This is not possible for multiple reasons: 
a)  $\Omega^1_X$ is not nef \cite[Corollary 2]{Yiyu}; b)  the Hurwitz formula gives $K_X=f^*(L)$, hence $X$ is of general type; c) up to a translation, a surjective map of tori of the same dimension is an isogeny, so in particular it is \'etale, while $f$ is ramified by construction. 
\end{itemize}
\end{pf}
 
%We conclude this section by pointing out that Theorem \ref{main} is a particular case of a much more general non-asphericity result that can be derived from the results of Pardini and Tovena \cite{PT95} and K\'ollar and Pardon \cite{Pardon}.

\begin{rem} The strategy of proof of  Theorem \ref{main} can be applied to the more general setting of a smooth abelian cover $f\colon X\to Y$ such that $Y$ is aspherical,  the  irreducible components  of the branch divisor are ample, and $f$ is completey ramified, i.e., it does not factor through a non-trivial \'etale cover $Y'\to Y$. Indeed, if this is the case, by \cite[Proposition 3.3]{PT95} the map $f_*\colon \pi_1(X)\to \pi_1(Y)$ is surjective with kernel a finite central subgroup of $\pi_1(Y)$, say $N$. If $N \ne \{1\}$, then $X$ is not aspherical for elementary topological reasons (see \cite[Proposition 2.45]{Hat}). If $N$ is trivial, then arguing with the degree of the map $f$ as in the proof of Theorem \ref{main} one can reach a contradiction.
\end{rem}

\begin{rem}\label{rem:  kollar}
After the appearance of this manuscript on the arXiv, J\'anos Koll\'ar pointed out to us that \cite[Theorem 16]{Pardon} can be used to prove the non-asphericity of the examples considered in Theorem \ref{main}. This approach is more general since it applies to certain analytic maps between compact complex analytic spaces where the codomain can be mildly singular. In our case one argues as follows:  $f_*\colon \pi_1(X)\to \pi_1(Y)$ is an isomorphism by \cite[Cor.~3.4]{PT95}, so if $X$ is aspherical we can  apply  \cite[Theorem 16]{Pardon} to  the map $f\colon X\to Y$ and deduce that it is a topological covering, contradicting the fact that  by construction $f$ is branched on a divisor of $Y$. Note that also one of the arguments given in the proof by rigidity  of Theorem \ref{main} leads to the same contradiction.
\end{rem}

\section{The Case of Surfaces}\label{surfaces}

Here we reprove Theorem \ref{main} in case $n=2$ by a different argument. The advantage of this approach is that it can be used to prove non-asphericity for many examples of surfaces of general type.
We start with a simple observation:

\begin{prop} \label{prop: pi2}
Let $X$ be a smooth $n$-dimensional complex projective variety and let $p\colon \widetilde X\to X$ be its universal cover. Let $B$ a smooth curve and $g\colon B\to X$   a non-constant morphism. 
Then:
\begin{enumerate}
\item if $g_*( \pi_1(B))\subset \pi_1(X)$ is a  finite subgroup, then $H_2(\widetilde X,\R)\ne 0$;
\item if $B=\mathbb P^1$, then $\pi_2(\widetilde X)$ is infinite.
\end{enumerate}
In particular, in both cases $X$ is not aspherical. 
\end{prop}
\begin{pf}
(1)  Let $H$ be an ample line bundle on $X$. Since $g$ is non constant we have $\deg_B(g^*H)>0$. If $\omega$ is a 2-form that represents the class of $H$ in  $H^2(X, \R)$, then $\deg_B(g^*H)=\int_Bg^*\omega$, so $B$ defines a non-zero integral class in $H_2( X,\R)$. Let now $B_0$ be a connected component of $B\times_X\widetilde X$. The induced map $q\colon B_0\to B$ is a finite \'etale morphism of degree equal to the cardinality  $\nu$ of $g_*( \pi_1(B))$ and $\int_{B_0}(g\circ q)^*\omega=deg_{B_0}((g\circ q)^*H)=\nu \deg_B(g^*H)>0$, so $B_0$ defines a non zero integral class in $H_2(\widetilde X,\R)$.
 \medskip
 
(2) By the proof of (1), the map $g$ gives a non-torsion  class in $H_2(X, \Z)$. By the Hurewicz homomorphism this implies that $g$ gives an element of infinite order of $\pi_2(X)$.
 \end{pf}

Let $\mathfrak M$ be a connected component of the moduli space of (canonical models of) surfaces of general type. It is well known that the minimal models of the surfaces parametrized by  $\mathfrak M$ are all diffeomorphic (see  \cite[Ch.~V]{tesi-manetti} for a nice discussion of this and other related facts).

% Indeed let $[Y]$ be a point of $\mathfrak M$ and let $\mathcal Y\to \Delta$ be a flat family over the complex disk such that for all $t\in \Delta$ the fiber $Y_t$ is the canonical model of a surface of general type and $Y_0\cong Y$. 
% The Brieskorn--Tyurina  simultaneous resolution  Theorem (cf. \cite[Theorem 3.1]{tesi-manetti}) \rita{qui \`e spiegato bene, gli articoli originali non li ho visti} tells us that, possibly after shrinking $\Delta$ and taking  a base change $\Delta\to \Delta$,  the family $\mathcal Y$ admits a simultaneous resolution. In other words, there is  a smooth family $\mathcal X\to\Delta$ and a birational morphism $\mathcal X \to \mathcal Y$ over $\Delta$ that restricts to the minimal resolution of the fiber $Y_t$ over each $t\in\Delta$. By Ehresmann's theorem all fibers $X_t$ are diffeomophic, so the diffeomorphism class of  the minimal model is locally constant, and therefore constant, on $\mathfrak M$.  
 
An immediate application of Proposition \ref{prop: pi2} is the following:
\begin{cor} \label{cor: M}
Let $\mathfrak M$ be a connected component of the moduli space  of canonical models of surfaces of general type.
 If a point of  $\mathfrak M$ corresponds to a singular surface, then  the minimal models of the surfaces  in $\mathfrak M$ have infinite $\pi_2$.
\end{cor}
\begin{pf} Let $[\bar S]\in \mathfrak M$ be a point corresponding to a singular surface  and let $S\to \bar S$ be the minimal resolution. Then $S$ has rational double points  and  the exceptional curves of $S\to \bar S$ are  $(-2)$-curves, namely they are  isomorphic to $\mathbb P^1$ and have self-intersection $-2$.  So $\pi_2(S)$ is infinite by Proposition \ref{prop: pi2}, (2). Since the minimal models of surfaces of $\mathfrak M$ are all diffeomorphic, this is enough to prove the statement.
\end{pf}

Variants of Corollary \ref{cor: M} can be used also when all surfaces in a connected component of the moduli are smooth. Indeed, arguing exactly as in the proof of Corollary \ref{cor: M}, one can show that if a surface $S$ with $K_S$ ample belongs a connected component  $\mathfrak M$ of the moduli space  that contains a surface $S'$ with a rational curve, then $S$ is not aspherical. We give an example below. 
 
\begin{ex} Consider the symmetric square  $S=S^2C$ with $C$ a curve of genus $g\ge 4$ (see  \cite[\S 2.4]{geography} for the definition and properties of symmetric squares).  The Albanese map of $S$ can be identified with the   map $a\colon S^2C\to\Pic^{(2)}(C)$ that associates to a degree 2 effective divisor on $C$ its linear equivalence class. So if $C$ is not hyperelliptic then the Albanese map is injective and therefore $S$ contains no rational curve. If $C$ is hyperelliptic  then the pairs of points in the canonical $g^1_2$ give a smooth rational curve $\Gamma$ on $S$, contracted by $a$, and $a$ is injective on $S\setminus \Gamma$. This implies that $\Gamma$ is the only rational curve contained in $S$. One can compute 
\[
K_S\Gamma=g-3\ge 1, \quad K^2_{S}=(g-1)(4g-9).
\]
Thus, Nakai's ampleness criterion (\cite[Corollary 6.4]{bpv})  implies that $K_S$ is ample in either case. Now, since all symmetric squares belong to the same irreducible component of the moduli space they are not aspherical. Interestingly, since 
\[
\chi(\OO_S)=\frac{g(g-3)}{2}+1,
\]
we have that the ratio $\frac{K^2_S}{\chi(\OO_S)}$ tends to $8$ from below as $g\to \infty$.
\end{ex}

Using this circle of ideas, and in particular Corollary \ref{cor: M}, we can now give an alternative proof of Theorem \ref{main} in dimension $n=2$.
\begin{prop} Let $S$ be a minimal surface of general type and maximal Albanese dimension  such that $K^2_S=4\chi(\OO_S)$. 
Then $S$ is not aspherical. \newline
In particular, the examples of Section  \ref{sec: main} are not aspherical for $n=2$.
\end{prop}
\begin{pf}
As explained in Remark \ref{rem: severi-line}, a surface $S$ as in the assumption is the minimal resolution of a double cover of an abelian surface $Y$ branched on an ample effective divisor $B$ with at most A-D-E singularities. If $B$ is smooth then $S$ is one of the examples constructed in Section  \ref{sec: main} while if $B$ is singular, then the double cover has canonical singularities and $S$ contains a rational curve. 
If $B$ is smooth, then by Corollary \ref{cor: M} it is enough to show that there is a singular surface in the  connected component of the moduli containing  $S$. Let the double cover be given by the relation $2L \sim B$ and let $D:=(d_1, d_2)$ be the type of the polarization $L$. The moduli space $\mathcal A_{2, D}$ of polarized abelian surfaces of type $(d_1, d_2)$  is irreducible (\cite[Thm.~8.2.6]{BL-CAV}) and it has a finite covering carrying a universal family (\cite[Thm.~7.9]{mumford-GIT}.  Using these facts it is a standard exercise to construct an irreducible  flat family of surfaces containing all the double covers of  abelian surfaces $Y$  with $L$ of type $D$ and $B$ with at most A-D-E singularities. This shows that  for fixed $D$ all our examples lie in the same component of the moduli space. 
This component  contains covers  where $Y=E_1\times E_2$ is   a product of elliptic curves,  $L$ is a product polarization and $B=B_1+B_2$, with $B_1$ the pull back of a smooth divisor of degree $2d_1$ on $E_1$ and $B_2$ the pullback  of a divisor of degree $2d_2$ on $E_2$. These  covers have  $A_1$ singularities at the intersection points of $B_1$ and $B_2$, so we can apply Corollary \ref{cor: M} and complete the proof.
\end{pf}

The case of surfaces is interesting for another reason. Indeed, Liu-Maxim-Wang conjecture implies a tantalizing conjecture regarding the geography of aspherical surfaces of general type. By the work of Demailly-Peternell-Schneider (see \cite[Theorem 2.5]{DPS94}) if $S$ is a K\"ahler surface with \emph{nef} holomorphic cotangent bundle, we then have that $c^2_1(S)\geq c_2(S)$. We therefore obtain the following statement.

\begin{prop}
Let $S$ be an aspherical surface of general type. If Liu-Maxim-Wang conjecture is true, we then have
\[
K^2_{S}\geq 6\chi(\mathcal{O}_S).
\]
\end{prop}

\begin{pf}
This follows by combining the standard identities
\[
K^2_{S}=c^2_{1}(S), \quad \chi(\mathcal{O}_S)=\frac{c^2_{1}(S)+c_2(S)}{12}
\]
with the Demailly-Peternell-Schneider inequality
\[
c^2_{1}(S)\geq c_{2}(S).
\]
\end{pf}

This proposition provides a concrete test on the plausibility of Liu-Maxim-Wang conjecture (see \cite[Conjecture 6.4]{LMW17b}). Indeed, any aspherical surface of general type $S$ with
\[
K^2_{S}< 6\chi(\mathcal{O}_S)
\]
would be a counterexample. Thus, it seems extremely interesting to ask if such surfaces exist. Unfortunately, the list of aspherical surfaces of general type seems to be not too rich. The list includes ball quotients, surfaces isogeneous to product of curves, Kodaira fibrations, Mostow-Siu surfaces, see the paper of Bauer-Catanese \cite{Fabrizio} for more details. In all of these examples, the stronger inequality
\begin{align}\label{bconjecture}
K^2_{S}\geq 8\chi(\mathcal{O}_S)
\end{align}
is satisfied. Maybe somewhat less well-known, the list of aspherical surfaces of general type includes also the vast majority of smooth minimal toroidal compactifications of ball quotients, see \cite[Theorem A]{tesi}. These examples tend to satisfy the bound in Equation \ref{bconjecture}, and it is indeed an interesting problem to determine whether Equation \ref{bconjecture} is satisfied by all aspherical toroidal compactifications of ball quotient surfaces. More generally, this circle of ideas highlights the search for aspherical surfaces of general type as an extremely interesting problem. Unfortunately, the vast majority of known  explicit constructions of surfaces of general type  involve linear systems and naturally specialize to examples with rational double points, and so cannot be aspherical by Corollary \ref{cor: M}. In the same vein, if a surface has a positive number of  moduli, also in this case it is reasonable to expect  that in most cases it can be specialized to a surface with rational double points and thus is not aspherical.  Since the expected number of moduli of a surface of general type, say $S$, is $10\chi(\OO_S)-2K_S$, it seems even harder to find examples of aspherical surfaces of general type such that 
\[
K^2_S/\chi(\OO_S)<5.
\]
Further exploration of these issues would seem to be one of the most compelling potential directions in the study of complex surfaces of general type and its interaction with low-dimensional topology. The many questions we highlighted here are a clear indication of the depth of our present ignorance and the frustrating lack of relevant examples.


\begin{thebibliography}{ELMNPM}
	%%%%%%%%%%%%%%%%%%%% 
	
\bibitem[AW21]{arapura} D. Arapura, B. Wang. Perverse sheaves on varieties with large fundamental group. \textit{arXiv:2109.07887v3[math.AG]}. To appear in \textit{J. Diff. Geom.}

\bibitem[Bar15]{barja-severi} M.A. Barja. Generalized Clifford-Severi inequality and the volume of irregular varieties. \textit{Duke Math. J. 164} (2015), no. 3, 541-568. 

\bibitem[BPS16]{BPS} M.A. Barja, R. Pardini, L. Stoppino. Surfaces on the Severi line. \textit{J. Math. Pures Appl. (9) 105} (2016), no. 5, 734-743. 

\bibitem[BHPV04]{bpv} W. P. Barth, K. Hulek, C. A. M. Peters, A. Van de Ven. Compact complex surfaces, 2nd ed., \textit{Ergebnisse der Mathematik und ihrer Grenzgebiete. 3. Folge}. A Series of Modern Surveys in Mathematics [Results in Mathematics and Related Areas. 3rd Series. A Series of Modern Surveys in Mathematics], vol. 4, Springer-Verlag, Berlin, 2004. 	

\bibitem[BC18]{Fabrizio} I. Bauer, F. Catanese. On rigid complex surfaces and manifolds  \textit{Adv. Math. 333}, (2018), 620-669. 

\bibitem[Ber03]{Berger} M. Berger. A panoramic view of Riemannian geometry. \textit{Springer-Verlag, Berlin},  2003. xxiv+824 pp.

\bibitem[BL04]{BL-CAV}  C. Birkenhake, H. Lange.  Complex abelian varieties. 2nd augmented ed. \textit{Grundlehren der Mathematischen Wissenschaften 302}. Berlin: Springer  xii, 635 p. (2004). 

\bibitem[Cat02]{Catanese} F. Catanese. Deformation types of real and complex manifolds, Contemporary trends in algebraic geometry and algebraic topology (Tianjin, 2000), 195-238, Nankai Tracts Math., 5, World Sci. Publishing, River Edge, NJ, 2002.

\bibitem[Che55]{Chern} S. S. Chern. On Curvature and Characteristic Classes of a Riemannian Manifold. \textit{Abh. Math. Semin. Univ. Hamb. 20} (1955), 117-126.

\bibitem[DPS94]{DPS94} J.- P. Demailly, T. Peternell, M. Schneider. Compact complex manifolds with numerically effective tangent bundles. \textit{J. Algebraic Geom. 3 } (1994), no. 2, 295–345.

\bibitem[DiC12]{tesi} L. F. Di Cerbo. Finite-volume complex-hyperbolic surfaces, their toroidal compactifications, and geometric applications. \textit{Pacific J. Math. 255} (2012), no. 2, 305-315.

\bibitem[DS22]{DS22} L. F. Di Cerbo, M. Stern. Price Inequalities and Betti numbers growth on manifolds without conjugate points. \textit{Comm. Anal. Geom. 30} (2022), no. 2, 297-334.

\bibitem[DX84]{DX84} H. Donnelly, F. Xavier. On the Differential Form Spectrum of Negatively Curved Riemannian Manifolds. \textit{Amer. J. Math. 106} (1984), no. 1, 169-185.

\bibitem [GL87]{GL1} M. Green, R. Lazarsfeld. Deformation theory, generic vanishing theorems and some conjectures of Enriques, Catanese and Beauville, \textit{Invent. Math. 90} (1987), 389--407.

\bibitem[Gro91]{Gro} M. Gromov. K\"ahler hyperbolicity and $L_2$-Hodge theory. \textit{J. Differential Geom. 33} (1991), no. 1, 263-292.

\bibitem[JZ00]{JZ} J. Jost, K. Zuo. Vanishing theorems for $L^2$-cohomology on infinite coverings of compact K\"ahler manifolds and applications in algebraic geometry. \textit{Comm. Anal. Geom. 8} (2000), no. 1, 1-30.

\bibitem[Hat02]{Hat} A. Hatcher. Algebraic Topology. \textit{Cambridge University Press, Cambridge}, 2002,  xii+544 pp.

\bibitem[Kol93]{Ko} J. Koll\'ar. Shafarevich maps and plurigenera of algebraic varieties. \textit{Invent. Math. 113} (1993), no. 1, 177-215.

\bibitem[Kol95]{KBook} J. Koll\'ar. Shafarevich maps and automorphic forms. M. B. Porter Lectures. \textit{Princeton University Press, Princeton, NJ}, 1995, x+201 pp.

\bibitem[KP12]{Pardon} J. Koll\'ar, J. Pardon. Algebraic varieties with semialgebraic universal cover. \textit{J. Topol. 5} (2012), no. 1, 199-212.

\bibitem[Laz04]{Laz} R. Lazarsfeld. Positivity in algebraic geometry I,textit{Ergebnisse der Mathematik und ihrer Grenzgebiete}, vol. 48, Berlin: Springer 2004.

\bibitem[LMW21]{LMW17b} Y. Liu, L. Maxim, B. Wang. Aspherical manifolds, Mellin transformation and a question of Bobadilla-Koll\'ar.  \textit{J. Reine Angew. Math. 781} (2021), 1-18.	 

\bibitem[Man96]{tesi-manetti} M. Manetti. Degeneration of algebraic surfaces and applications to moduli problems, \textit{edizioni Scuola Normale Superiore (EN) collana Tesi}, 1996. (available at https://www1.mat.uniroma1.it/people/manetti/dispense/tesiperf.pdf)

\bibitem[Max23]{Maxim} L. Maxim. On singular generalizations of the Singer–Hopf conjecture. \textit{Math. Nachr. 296} (2023), no. 11, 5232-5241.

\bibitem[MLP12]{geography} M. Mendes Lopes, R. Pardini, The geography of irregular surfaces, \textit{Current developments in algebraic geometry,  Math. Sci. Res. Inst. Publ. 59}, Cambridge Univ. Press (2012), 349--378. 

\bibitem[MFK94]{mumford-GIT}  D. Mumford, J.Fogarty, F. Kirwan, Geometric invariant theory. \textit{3rd enl. ed. Ergebnisse der Mathematik und ihrer Grenzgebiete. 2. Folge. 34.} Berlin: Springer-Verlag. 320 p. (1994). 

\bibitem[PT95]{PT95} R. Pardini, F. Tovena. On the fundamental group of an abelian cover. \textit{Internat. J. Math 6} (1995), no. 5, 767-789.

\bibitem[PS13]{PS13} M. Popa, Ch. Schnell, Generic vanishing theory via mixed Hodge modules. \textit{Forum Math. Sigma 1} (2013), e1, 60pp.

\bibitem[SY94]{SchoenY} R. Schoen, S.-T. Yau.  Lectures on differential geometry. Conference Proceedings and Lecture Notes in Geometry and Topology, I. \textit{International Press, Cambridge, MA}, 1994.

\bibitem[Wan22]{Yiyu} Y. Wang. Ramified covers of varieties with nef cotangent bundle. \textit{C. R. Math. Acad. Paris 360} (2022), 929-932.

\bibitem[Zha14]{zhang-severi} T. Zhang. Severi inequality for varieties of maximal Albanese dimension.  \textit{Math. Ann. 359} (2014), no. 3-4, 1097-1114.
\end{thebibliography}
\end{document}